\documentclass[reqno,a4paper]{amsart}

\usepackage{amsmath}
\usepackage{amssymb}
\usepackage{amsthm}
\usepackage{enumerate}
\usepackage{mathrsfs} 
\usepackage{eqlist}
\usepackage{array}
\usepackage{ifpdf}
\ifpdf
 \usepackage[hyperindex]{hyperref}
\else
 \expandafter\ifx\csname dvipdfm\endcsname\relax
 \usepackage[hypertex,hyperindex]{hyperref}%
 \else
 \usepackage[dvipdfm,hyperindex]{hyperref}%
 \fi
\fi

\theoremstyle{plain}
\newtheorem{theorem}{Theorem}[section]

\DeclareMathOperator{\td}{d}
\DeclareMathOperator{\ti}{i}
\DeclareMathOperator{\te}{e}
\DeclareMathOperator{\bell}{B}
\DeclareMathOperator{\TT}{T}

\theoremstyle{definition}

\newtheorem{remark}{\textup{Remark}} 

\numberwithin{equation}{section}

\begin{document}
	
	\title[Parametric kinds of generalized Apostol--Bernoulli polynomials]%
	{Parametric kinds of generalized Apostol--Bernoulli polynomials and their properties}
	\author[Zeynep \"Ozat, Bayram \c{C}ekim, Can K{\i}z{\i}late\c{s} \and Feng Qi]%
	{Zeynep \"Ozat$^{1}$, Bayram \c{C}ekim$^{1}$*, Can K{\i}z{\i}late\c{s}$^{2}$ \and Feng Qi$^{3}$}
	
	\newcommand{\acr}{\newline\indent}
	
	\address{\llap{$^{1}$\,}Department of Mathematics \acr
			Gazi University \acr
			Teknikokullar, 06500, Ankara, Turkey}
	\email{zeynepozat95@gmail.com, bayramcekim@gazi.edu.tr}
	
	\address{\llap{$^{2}$\,}Department of Mathematics \acr
			 Zonguldak Bülent Ecevit University \acr
			  67100, Zonguldak, Turkey}
	\email{cankizilates@gmail.com}
	
	\address{\llap{$^{3}$\,}Institute of Mathematics\acr
		 Henan Polytechnic University \acr
		  Jiaozuo 454010, Henan, China \acr
	  Independent Researcher, Dallas, TX 75252-8024, USA}
	\email{qifeng618@gmail.com, qifeng618@hotmail.com, qifeng618@qq.com}
	
	\thanks{*Corresponding author}
	
	
	\begin{abstract}
	The purpose of this paper is to define generalized Apostol--Bernoulli polynomials with including a new cosine and sine parametric type of generating function using the quasi-monomiality properties and trigonometric functions. In this study, the Apostol-Bernoulli polynomials with three variable  are defined with two new generating functions cosine and sine parameters. Then, we investigate multiplicative and derivative operators, diffrential equations, some summation formulas and partial differential equations for these polynomials. Moreover, we introduce Gould--Hopper--Apostol--Bernoulli type polynomials, Hermite--Appell--Apostol--Bernoulli type polynomials and truncated exponential Apostol--Bernoulli type polynomials. Finally, the special cases of these new polynomials are investigated, and the corresponding results are expressed.	
	\end{abstract}
	
	\maketitle
	\section{Introduction}
	Generating functions play an important role in the investigation of properties of special polynomials and numbers. There are many studies related to polynomials and their generating functions. A few of references to special polynomials and their generating functions are given in the monographs~\cite{atf3, atf4, atf2, atf1}. On the other hand, generating functions have some applications in many fields such as applied mathematics, algebra, statistics, combinatorics, and physics.
	\par
	The Bernoulli polynomials are one of the important study subject in mathematics. These polynomials are generated by using exponential functions and series expansions. These polynomials have many application areas such as applied mathematics, computer programming, and statistics. The Bernoulli numbers and polynomials are related to many special numbers and polynomials such as the Euler, Genocchi, Hermite, and Stirling polynomials and numbers. The Swiss mathematician Jacob Bernoulli first mentioned the Bernoulli numbers in his work~\cite{atf5}. Then, Tom M. Apostol developed Bernoulli polynomials in~\cite{tm1}. Later, some properties and relations related to these polynomials are obtained~\cite{b5, fq14, fq5, fq4, fq3, fq11, fq2, fq7, fq, fq1, fq16, fq10, Nice-ApBern-Polyn-2v.tex, fq8, fq15, fq6, fq17, fq13, JNSA201701222.tex}.
	\par
	Throughout this paper, we use the notions $\mathbb{N}=\{1,2,3,\dotsc\}$, $\mathbb{N}_0=\{0\}\cup\mathbb{N}$, $\mathbb{R}$, and $\mathbb{C}$ to denote the sets of positive integers numbers, non-negative integers numbers, real numbers, and complex numbers, respectively.
	\par
	The classical Bernoulli polynomials $\bell_{n}(x)$ were defined in~\cite[p.~1]{b5} and~\cite[p.~3]{atf3} by the generating function
	\begin{equation*}
		\sum_{n=0}^\infty\bell_{n}(x)\frac{t^{n}}{n!}
		=\frac{t\te^{xt}}{\te^{t}-1},\quad |t|<2\pi.
	\end{equation*}
	The first five Bernoulli polynomials are
	\begin{gather*}
		\bell_0(x)=1,\quad \bell_1(x)=x-\frac12,\quad \bell_2(x)=x^2-x+\frac16,\\
		\bell_3(x)=x^3-\frac32x^2+\frac12x,\quad \bell_4(x)=x^4-2x^3+x^2-\frac1{30}.
	\end{gather*}
	\par
	The quantities $\bell_{n}(0)=\bell_n$ are called the Bernoulli numbers.
	Due to that the function $\frac{t}{\te^t-1}-1+\frac{t}2$ is even in $t\in\mathbb{R}$, all of the Bernoulli numbers $\bell_{2n+1}$ for $n\in\mathbb{N}$ equal $0$. The first few Bernoulli numbers $\bell_{2n}$ are
	\begin{gather*}
		\bell_2=\frac16,\quad \bell_4=-\frac1{30},\quad \bell_6=\frac1{42},\quad \bell_8=-\frac1{30},\quad \bell_{10}=\frac5{66},\\ \bell_{12}=-\frac{691}{2730},\quad \bell_{14}=\frac76,\quad \bell_{16}=-\frac{3617}{510},\quad \bell_{18}=\frac{43867}{798}.
	\end{gather*}
	In~\cite{fq2, fq13}, some new properties for the ratio $\frac{\bell_{2n+2}}{\bell_{2n}}$ with $n\in\mathbb{N}$ were discovered.
	\par
	One of the relations between the Bernoulli numbers and polynomials is
	$$
	\bell_{n}(x)=\sum_{k=0}^n\binom{n}{k}\bell_kx^{n-k},\quad n\in\mathbb{N}_0
	$$
	which is a special case of the identity
	\begin{equation*}
		\bell_{n}(x+h)=\sum_{k=0}^n\binom{n}{k}\bell_k(x)h^{n-k},\quad n\in\mathbb{N}_0
	\end{equation*}
	in~\cite[p.~802, 23.1.7]{abram}.
	\par
	The generalized Bernoulli polynomials $\bell_{n}^{(v)}(x)$ of order $v$ are defined~\cite[p.~4]{Temme-96-book} by
	\begin{equation*}
		\sum_{n=0}^\infty\bell_{n}^{(v)}(x) \frac{t^{n}}{n!}=\biggl(\frac{t}{\te^{t}-1}\biggr)^v \te^{xt},\quad |t|<2\pi.
	\end{equation*}
	It is easy to see that $\bell_{n}^{(1)}(x)=\bell_{n}(x)$. In addition, we call $\bell_{n}^{(v)}=\bell_{n}^{(v)}(0)$ the Bernoulli numbers of order $v$. The relation between $\bell_{n}^{(v)}$ and $\bell_{n}^{(v)}(x)$ was given~\cite[p.~510]{brn} by
	\begin{equation*}
		\bell_{n}^{(v)}(x)=\sum_{k=0}^n\binom{n}{k}\bell_k^{(v)}x^{n-k}.
	\end{equation*}
	\par
	The Apostol--Bernoulli polynomials $\bell_{n}(x;\lambda)$ were defined in~\cite[p.~165]{tm1} by
	\begin{equation*}
		\sum_{n=0}^\infty\bell_{n}(x;\lambda)\frac{t^{n}}{n!}=\frac{t\te^{xt}}{\lambda\te^{t}-1},
	\end{equation*}
	where $|t|<2\pi$ when $\lambda=1$ or $|t|<|\ln\lambda|$ when $\lambda\ne 1$ and $\lambda\in\mathbb{C}$.
	When $x=0$, we call $\bell_{n}(0;\lambda)=\bell_{n}(\lambda)$ the Apostol--Bernoulli numbers.
	These polynomials satisfy
	\begin{align*}
		\bell_{n}(x;\lambda)&=\sum_{k=0}^n\binom{n}{k}\bell_k(\lambda)x^{n-k},\\
		\bell_{n}(x+y;\lambda)&=\sum_{k=0}^n\binom{n}{k}\bell_k(x;\lambda)y^{n-k},\\
		\frac{\partial^p}{\partial x^p}\bell_{n}(x;\lambda)&=\frac{n!}{(n-p)!}\bell_{n-p}(x;\lambda),\\
		\int_{a}^{b}\bell_{n}(t;\lambda)\td t&=\frac{\bell_{n+1}(b;\lambda)-\bell_{n+1}(a;\lambda)}{n+1},\\
		\sum_{k=0}^{m-1}k^n&=\frac{\lambda-1}{n+1}\sum_{k=1}^m\bell_{n+1}(k;\lambda)+\frac{\bell_{n+1}(m;\lambda)-\bell_{n+1}(\lambda)}{n+1}
	\end{align*}
	for $m,n\in\mathbb{N}_0$ and $0\le p\le n$.
	In~\cite[p.~165]{tm1}, Apostol gave the relation
	$$
	\lambda\bell_n(x+1,\lambda)-\bell_n(x,\lambda)=nx^{n-1},\quad n\ge1.
	$$
	The first six Apostol--Bernoulli numbers are
	\begin{gather*}
		\bell_0(\lambda)=0,\quad
		\bell_1(\lambda)=\frac{1}{\lambda-1},\quad
		\bell_2(\lambda)=-\frac{2\lambda}{(\lambda-1)^2},\quad
		\bell_3(\lambda)=\frac{3\lambda(\lambda+1)}{(\lambda-1)^3},\\
		\bell_4(\lambda)=\frac{4\lambda(\lambda^{2}+4\lambda+1)}{(\lambda-1)^4},\quad
		\bell_5(\lambda)=\frac{5\lambda(\lambda^{3}+11\lambda^{2}+11\lambda+1)}{(\lambda-1)^5}.
	\end{gather*}
	\par
	The generalized Apostol--Bernoulli polynomials $\bell_{n}^{(v)}(x;\lambda)$ of order $v$ can be generated by
	\begin{equation}\label{a-p}
		\sum_{n=0}^\infty\bell_{n}^{(v)}(x;\lambda)\frac{t^{n}}{n!}
		=\biggl(\frac{t}{\lambda\te^{t}-1}\biggr)^{v}\te^{xt},\quad |t+\ln\lambda|<2\pi
	\end{equation}
	with
	$$
	\bell_{n}^{(v)}(x)=\bell_{n}^{(v)}(x;1) \quad \text{and} \quad \bell_{n}^{(v)}(\lambda)=\bell_{n}^{(v)}(0;\lambda),
	$$
	where $\bell_{n}^{(v)}(\lambda)$ denote the Apostol--Bernoulli numbers of order $v$, see~\cite[p.~292]{atf20} and~\cite[p.~92]{atf2}. The identity between the Apostol--Bernoulli polynomials $\bell_{n}^{(v)}(x;\lambda)$ of order $v$ and the Apostol--Bernoulli numbers $\bell_{n}^{(v-1)}(\lambda)$ of order $v-1$ is given in~\cite[p.~300]{atf20} by
	$$
	\bell_{n}^{(v)}(x;\lambda)=\sum_{k=0}^n\binom{n}{k}\bell_{n-k}^{ (v-1)}(\lambda)\bell_{k}(x;\lambda),\quad n\in\mathbb{N}_0.
	$$
	The Apostol--Bernoulli polynomials $\bell_{n}^{(v)}(x;\lambda)$ of order $v$ satisfy
	\begin{align*}
		\bell_{n}^{(0)}(x;\lambda)&=x^n,\\
		\bell_{n}^{(v)}(x;\lambda)&=\sum_{k=0}^n\binom{n}{k}\bell_{k}^{(v)}(\lambda)x^{n-k},\\
		\lambda\bell_n^{(v)}(x+1;\lambda)-\bell_n^{(v)}(x;\lambda)&=n\bell_{n-1}^{(v-1)}(x;\lambda),\\
		\bell_{n}^{(v+u)}(x+y;\lambda)&=\sum_{k=0}^n\binom{n}{k}\bell_k^{(v)}(x;\lambda)\bell_{n-k}^{(u)}(y;\lambda),\\
		\frac{\partial}{\partial x}\bell_{n}^{(v)}(x;\lambda)&=n\bell_{n-1}^{(v)}(x;\lambda),\\
		\int_{a}^{b}\bell_{n}^{(v)}(t;\lambda)\td t&=\frac{\bell_{n+1}^{(v)}(b;\lambda)-\bell_{n+1}^{(v)}(a;\lambda)}{n+1}.
	\end{align*}
	The properties of Apostol--Bernoulli polynomials $\bell_{n}^{(v)}(x;\lambda)$ of order $v$ are not limited to these. In the literature, there are many studies on the Apostol--Bernoulli polynomials such as Gaussian hypergeometric function, multiplication formulas, and fourier expansions. The relations between the Apostol--Bernoulli polynomials, the Euler polynomials, the Genocchi polynomials, and the Stirling numbers have also been investigated in~\cite{ab6, ab9, ab0, ab4, atf20, lu1, ab7}.
	\par
	In~\cite[p.~3]{Jamei-Beyki-Koepf-MJM2018} and~\cite[p.~4]{ms1}, the Maclaurin expansions of two functions $\te^{xt}\cos(yt)$ and $\te^{xt}\sin(yt)$ are given by
	\begin{equation}\label{1.5}
		\te^{xt}\cos(yt)=\sum_{n=0}^\infty C_{n}(x,y)\frac{t^{n}}{n!}
	\end{equation}
	and
	\begin{equation}\label{1.6}
		\te^{xt}\sin(yt)=\sum_{n=0}^\infty S_{n}(x,y)\frac{t^{n}}{n!}
	\end{equation}
	respectively, where
	\begin{align*}
		C_{n}(x,y)&=\sum_{r=0}^{\lfloor n/2\rfloor}(-1)^{r}\binom{n}{2r}x^{n-2r}y^{2r},\\
		S_{n}(x,y)&=\sum_{r=0}^{\lfloor(n-1)/2\rfloor}(-1)^{r}\binom{n}{2r+1}x^{n-2r-1}y^{2r+1},
	\end{align*}
	and $\lfloor\lambda\rfloor$ denotes the floor function defined by the largest integer less than or equal to $\lambda\in\mathbb{R}$.
	\par
	In~\cite[p.~909]{ms2}, Srivastava et al defined two parametric types Apostol--Bernoulli polynomials $\bell_{n}^{(c,v)}(x,y;\lambda)$ and $\bell_{n}^{(s,v)}(x,y;\lambda)$ by
	\begin{equation}\label{1.7}
		\biggl(\frac{t}{\lambda\te^{t}-1}\biggr)^{v}\te^{xt}\cos(yt)=\sum_{n=0}^\infty\bell_{n}^{(c,v)}(x,y;\lambda)\frac{t^{n}}{n!}
	\end{equation}
	and
	\begin{equation}\label{1.7.7}
		\biggl(\frac{t}{\lambda\te^{t}-1}\biggr)^{v}\te^{xt}\sin(yt)=\sum_{n=0}^\infty\bell_{n}^{(s,v)}(x,y;\lambda)\frac{t^{n}}{n!}.
	\end{equation}
	\par
	For $v\in\mathbb{N}_{0}$ and $|t|<\ln2$, Srivastava and K{\i}z{\i}late\c{s} introduced in~\cite[p.~3257]{kz1} two parametric kinds of the Fubini type polynomials $F_{n}^{(c,v)} (x,y)$ and $F_{n}^{(s,v)} (x,y)$ by
	$$
	\frac{2^{v}}{(2-\te^{t})^{2v}}\te^{xt} \cos(yt)=\sum_{n=0}^\infty F_{n}^{(c,v)} (x,y) \frac{t^{n}}{n!},
	$$
	and
	$$
	\frac{2^{v}}{(2-\te^{t})^{2v}}\te^{xt}\sin(yt)=\sum_{n=0}^\infty F_{n}^{(s,v)}(x,y) \frac{t^{n}}{n!}.
	$$
	\par
	Some general polynomials have been described with the help of monomiality principle. Khan and Raza defined in~\cite[p.~4]{rz1} the $2$-variable general polynomials $\TT_{n}(x,y)$ by
	\begin{equation}\label{1.9}
		\te^{xt}\mathcal{U}(y,t)=\sum_{n=0}^\infty\TT_{n}(x,y)\frac{t^{n}}{n!},\quad\TT_{0}(x,y)=1,
	\end{equation}
	where
	\begin{equation}\label{1.10}
		\mathcal{U}(y,t)=\sum_{n=0}^\infty\mathcal{U}_{n}(y)\frac{t^{n}}{n!},\quad\mathcal{U}_{0}(y)\ne0.
	\end{equation}
	The multiplicative and derivative operators for the polynomials $\TT_{n}(x,y)$ are defined in~\cite[p.~4]{rz1} by
	\begin{equation}\label{1.11}
		\hat{M}_{\TT}=x+\frac{\mathcal{U}'(y,D_{x})}{\mathcal{U}(y,D_{x})}
		\quad\text{and}\quad
		\hat{P}_{\TT}=D_{x}.
	\end{equation}
	According to monomiality principle, the polynomials $\TT_{n}(x,y)$ satisfy
	\begin{align}
		\hat{M}_{\TT}\{\TT_{n}(x,y)\}&=\TT_{n+1}(x,y),\label{1.13}\\
		\hat{P}_{\TT}\{\TT_{n}(x,y)\}&=n\TT_{n-1}(x,y),\label{1.14}\\
		\hat{M}_{\TT}\hat{P}_{\TT}\{\TT_{n}(x,y)\}&=n\TT_{n}(x,y),\label{1.15}\\
		\exp\bigl(\hat{M}_{\TT}t\bigr)&=\sum_{n=0}^\infty\TT_{n}(x,y)\frac{t^{n}}{n!},\quad|t|<\infty.\label{1.16}
	\end{align}
	\par
	Among other things, via monomiality principle, Srivastava and five coauthors defined in~\cite[p.~5]{sr1} two parametric types of the generalized Fubini type polynomials by
	\begin{equation*}
		\frac{2^{v}}{(2-\te^{t})^{2v}}\te^{xt}\mathcal{U}(y,t)\cos(zt)=\sum_{n=0}^\infty\mathcal{F}_{n}^{(c,v)}(x,y,z)\frac{t^{n}}{n!}
	\end{equation*}
	and
	\begin{equation*}
		\frac{2^{v}}{(2-\te^{t})^{2v}}\te^{xt}\mathcal{U}(y,t)\sin(zt)=\sum_{n=0}^\infty\mathcal{F}_{n}^{(s,v)}(x,y,z)\frac{t^{n}}{n!},
	\end{equation*}
	where $v\in\mathbb{N}_{0}$ and $|t|<\ln2$.
	\par
	This paper is organized as follows. In Section~\ref{section2}, by using operational methods and monomiality principle, we define the parametric kinds of generalized Apostol--Bernoulli polynomials and obtain quasi-monomial properties for these kinds of polynomials. In Section~\ref{section3}, we establish various relations, identities, and summation formulas for generalized Apostol--Bernoulli polynomials. In Section~\ref{section4}, we discover partial differential equations and identities for their generating functions. In Section~\ref{section5}, using the generating functions for the parametric kinds of generalized Apostol--Bernoulli polynomials, we present new polynomials and investigate their properties.
	
	\section{Parametric kinds of generalized Apostol--Bernoulli polynomials}\label{section2}
	
	In this section, by using monomiality principle and operational methods, via generating functions, we give two parametric types of generalized Apostol--Bernoulli type polynomials and introduce multiplicative and derivative operators and differential equations.
	
	\begin{theorem}\label{them-2.1-2.2}
		The generating functions for the parametric kinds of generalized Apostol--Bernoulli polynomials ${}_{\TT\!}\bell_{n}^{(c,v)}(x,y,z;\lambda)$ and ${}_{\TT\!}\bell_{n}^{(s,v)}(x,y,z;\lambda)$ are defined
		\begin{equation}\label{2.1}
			\biggl(\frac{t}{\lambda\te^{t}-1}\biggr)^{v}\te^{xt}\mathcal{U}(y,t) \cos(zt)
			=\sum_{n=0}^\infty{}_{\TT\!}\bell_{n}^{(c,v)}(x,y,z;\lambda) \frac{t^{n}}{n!},
		\end{equation}
		and
		\begin{equation}\label{2.2}
			\biggl(\frac{t}{\lambda\te^{t}-1}\biggr)^{v}\te^{xt}\mathcal{U}(y,t) \sin(zt)
			=\sum_{n=0}^\infty{}_{\TT\!}\bell_{n}^{(s,v)}(x,y,z;\lambda) \frac{t^{n}}{n!},
		\end{equation}
		where
		$$ \vert t \vert <2\pi, \,\, \text{when}  \,\, \lambda=1 ; \, \vert t \vert < \vert \log \lambda \vert,  \,\, \text{when}  \,\, \lambda \neq 1; \,\,  1^{v}:=1. $$
	\end{theorem}
	
	\begin{proof}
		Replacing $x$ and $y$ by the multiplication operator $\hat{M}_{\TT}$ and $z$ in~\eqref{a-p} gives
		\begin{equation*}
			\biggl(\frac{t}{\lambda\te^{t}-1}\biggr)^{v}\exp\bigl(\hat{M}_{\TT}t\bigr)\mathcal{U}(y,t)\cos(zt)
			=\sum_{n=0}^\infty\bell_{n}^{(c,v)}(x,y,z;\lambda)\frac{t^{n}}{n!}.
		\end{equation*}
		Using equalities in~\eqref{1.11} and~\eqref{1.16} results in
		\begin{equation*}
			\biggl(\frac{t}{\lambda\te^{t}-1}\biggr)^{v}\Biggl[\sum_{n=0}^\infty\TT_{n}(x,y)\frac{t^{n}}{n!}\Biggr]\cos(zt)
			=\sum_{n=0}^\infty{}\bell_{n}^{(c,v)}\biggl(x+\frac{\mathcal{U}'(y,D_x)}{\mathcal{U}(y,D_x)},z;\lambda\biggr)\frac{t^{n}}{n!}.
		\end{equation*}
		Using the equation~\eqref{2.1} and denoting the resultant parametric kinds of generalized Apostol--Bernoulli polynomials by ${}_{\TT\!}\bell_{n}^{(c,v)}(x,y,z;\lambda)$ yield
		\begin{equation*}
			{}_{\TT\!}\bell_{n}^{(c,v)}(x,y,z;\lambda)
			=\bell_{n}^{(c,v)}\biggl(x+\frac{\mathcal{U}'(y,D_x)}{\mathcal{U}(y,D_x)},z;\lambda\biggr).
		\end{equation*}
		Therefore, we acquire the assertion in the equation~\eqref{2.1}.
		\par
		Similarly, we can prove the assertion in the equation~\eqref{2.2}.
		The proof of Theorem~\ref{them-2.1-2.2} is complete.
	\end{proof}
	
	In order to find the quasi-monomial properties of parametric kinds of generalized Apostol--Bernoulli polynomials ${}_{\TT\!}\bell_{n}^{(c,v)}(x,y,z;\lambda)$ and ${}_{\TT\!}\bell_{n}^{(s,v)}(x,y,z;\lambda)$, we prove the following properties.
	
	\begin{theorem}\label{thm-2.3-2.6}
		The parametric kinds of generalized Apostol--Bernoulli polynomials ${}_{\TT}\!\bell_{n}^{(c,v)}(x,y,z;\lambda)$ and ${}_{\TT\!}\bell_{n}^{(s,v)}(x,y,z;\lambda)$ are quasi-monomial with respect to the following multiplicative and derivative operators:
		\begin{align}
			\hat{M}_{\TT\bell c}&=x+\frac{\mathcal{U}'(y,D_x)}{\mathcal{U}(y,D_x)} +\frac{v}{D_{x}}\frac{\lambda\te^{D_{x}}(1-D_{x})-1} {\lambda\te^{D_{x}}-1}-z\tan(zD_{x})\label{2.3},\\
			\hat{P}_{\TT\bell c}&=D_{x},\label{2.4}\\
			\hat{M}_{\TT\bell s}&=x+\frac{\mathcal{U}'(y,D_x)}{\mathcal{U}(y,D_x)}+\frac{v}{D_{x}} \frac{\lambda\te^{D_{x}}(1-D_{x})-1}{\lambda\te^{D_{x}}-1}+z\cot(zD_{x}),\label{2.5}\\
			\hat{P}_{\TT\bell s}&=D_{x}.\label{2.6}
		\end{align}
	\end{theorem}
	
	\begin{proof}
		Differentiating the equation~\eqref{2.1} with respect to $t$ leads to
		\begin{equation*}
			\biggl[x+\frac{\mathcal{U}'(y,t)}{\mathcal{U}(y,t)}+\frac{v}{t}\frac{\lambda\te^{t}(1-t)-1}{\lambda\te^{t}-1}\biggr] \biggl(\frac{t}{\lambda\te^{t}-1}\biggr)^{v}\te^{xt}\mathcal{U}(y,t) \cos(zt)
			=\sum_{n=0}^\infty{}_{\TT\!}\bell_{n+1}^{(c,v)}(x,y,z;\lambda) \frac{t^{n}}{n!},
		\end{equation*}
		which, in view of the equation~\eqref{2.1}, is equivalent to
		\begin{equation}\label{2.7}
			\sum_{n=0}^\infty\biggl[x+\frac{\mathcal{U}'(y,t)}{\mathcal{U}(y,t)}+\frac{v}{t}\frac{\lambda\te^{t}(1-t) -1}{\lambda\te^{t}-1}\biggr]\biggl[{}_{\TT\!}\bell_{n}^{(c,v)}(x,y,z;\lambda) \frac{t^{n}}{n!}\biggr]
			=\sum_{n=0}^\infty{}_{\TT\!}\bell_{n+1}^{(c,v)}( x,y,z;\lambda)\frac{t^{n}}{n!}.
		\end{equation}
		If $\mathcal{U}(y,t)$ is invertible and $\frac{\mathcal{U}'(y,t)}{\mathcal{U}(y,t)}$ has Taylor's series expansion in powers of $\frac{t^{n}}{n!}$, applying the identity~\eqref{2.7}, we obtain
		\begin{equation*}
			D_{x}\bigl[\te^{xt}\mathcal{U}(y,t)\bigr]=t\bigl[\te^{xt}\mathcal{U}(y,t)\bigr].
		\end{equation*}
		In the equation~\eqref{2.7}, comparing the coefficients of $t^n$ on both sides results in
		\begin{equation*}
			\biggl[x+\frac{\mathcal{U}'(y,D_x)}{\mathcal{U}(y,D_x)}+\frac{v}{D_{x}}\frac{\lambda\te^{D_{x}}(1-D_{x})-1} {\lambda\te^{D_{x}}-1}-z\tan(zD_{x})\biggr] \bigl[{}_{\TT}\!\bell_{n}^{(c,v)}(x,y,z;\lambda)\bigr]
			={}_{\TT\!}\bell_{n+1}^{(c,v)}(x,y,z;\lambda),
		\end{equation*}
		which, in view of the monomiality principle given in the equation~\eqref{1.13}, means the assertion in the equation~\eqref{2.3}.
		\par
		Differentiating on both sides of the equation~\eqref{2.1} with respect to $x$ gives
		\begin{equation*}
			D_{x}\Biggl[\sum_{n=0}^\infty{}_{\TT}\!\bell_{n}^{(c,v)}(x,y,z;\lambda) \frac{t^{n}}{n!}\Biggr]
			=\sum_{n=0}^\infty{}_{\TT}\!\bell_{n-1}^{(c,v)}(x,y,z;\lambda) \frac{t^{n}}{(n-1) !}.
		\end{equation*}
		When comparing the coefficients of $\frac{t^n}{n!} $ on both sides, we obtain the monomiality principle property given in the equation~\eqref{1.14} for the equation~\eqref{2.4}.
		\par
		By virtue of the method as in the proof of the equations~\eqref{2.3} and~\eqref{2.4}, we can prove the claims in the equations~\eqref{2.5} and~\eqref{2.6}. The proof of Theorem~\ref{thm-2.3-2.6} is complete.
	\end{proof}
	
	\begin{theorem}\label{thm-2.8-2.9}
		The parametric kinds of generalized Apostol--Bernoulli polynomials ${}_{\TT\!}\bell_{n}^{(c,v)}(x,y,z;\lambda)$ and ${}_{\TT\!}\bell_{n}^{(s,v)}(x,y,z;\lambda)$ satisfy
		\begin{equation}\label{2.8}
			\biggl[xD_{x}+\frac{\mathcal{U}'(y,D_x)}{\mathcal{U}(y,D_x)}D_{x}+v\frac{\lambda\te^{D_{x}}(1-D_{x}) -1
			}{\lambda\te^{D_{x}}-1}-z\tan (zD_{x}) D_{x}-n\biggr] {}_{\TT}\!\bell_{n}^{(c,v)}(x,y,z;\lambda)=0
		\end{equation}
		and
		\begin{equation}\label{2.9}
			\biggl[xD_{x}+\frac{\mathcal{U}'(y,D_x)}{\mathcal{U}(y,D_x)}D_{x}+v\frac{\lambda\te^{D_{x}}(1-D_{x}) -1}{\lambda\te^{D_{x}}-1}+z\cot (zD_{x}) D_{x}-n\biggr] {}_{\TT}\!\bell_{n}^{(s,v)}(x,y,z;\lambda)=0.
		\end{equation}
	\end{theorem}
	
	\begin{proof}
		Using operators in the equations~\eqref{2.3} and~\eqref{2.4} and employing the monomiality principle in the equation~\eqref{1.15}, we obtain the assertion in the equation~\eqref{2.8}. Similarly, we can prove the result in the equation~\eqref{2.9}. The proof of Theorem~\ref{thm-2.8-2.9} is complete.
	\end{proof}
	
	\section{Identities and relations}\label{section3}
	In this section, by using generating functions, we derive some identities and recursive relations for the parametric kinds of generalized Apostol--Bernoulli polynomials.
	
	\begin{theorem}\label{thm-3.1-3.2}
		The parametric kinds of generalized Apostol--Bernoulli polynomials ${}_{\TT\!}\bell_{n}^{(c,v)}(x,y,z;\lambda)$ and ${}_{\TT\!}\bell_{n}^{(s,v)}(x,y,z;\lambda)$ satisfy
		\begin{equation}{\label{3.1}}
			{}_{\TT\!}\bell_{n}^{(c,v)}(x,y,z;\lambda)=\sum_{r=0}^n\binom{n}{r}\bell_{n-r}^{(c,v)}(0,z;\lambda) \TT_{r}(x,y)
		\end{equation}
		and
		\begin{equation}{\label{3.2}}
			{}_{\TT\!}\bell_{n}^{(s,v)}(x,y,z;\lambda)=\sum_{r=0}^n\binom{n}{r}\bell_{n-r}^{(s,v)}(0,z;\lambda) \TT_{r}(x,y).
		\end{equation}
	\end{theorem}
	
	\begin{proof}
		From~\eqref{1.9} and~\eqref{1.7}, we derive
		\begin{equation*}
			\sum_{n=0}^\infty{}_{\TT\!}\bell_{n}^{(c,v)}(x,y,z;\lambda) \frac{t^{n}}{n!}
			=\sum_{n=0}^\infty\sum_{r=0}^n \binom{n}{r}\bell_{n}^{(c,v)}(0,z;\lambda) \TT_{r}(x,y)\frac{t^{n}}{n!}.
		\end{equation*}
		Equating coefficients of $\frac{t^n}{n!}$ on both sides of this equation yields the assertion in the equation~\eqref{3.1}.
		Similarly, we can verify the assertion in the equation~\eqref{3.2}.
		The proof of Theorem~\ref{thm-3.1-3.2} is complete.
	\end{proof}
	
	\begin{theorem}\label{thm3.3-3.4}
		The parametric kinds of generalized Apostol--Bernoulli polynomials ${}_{\TT\!}\bell_{n}^{(c,v)}(x,y,z;\lambda)$ and ${}_{\TT\!}\bell_{n}^{(s,v)}(x,y,z;\lambda)$ satisfy
		\begin{equation}\label{3.3}
			{}_{\TT\!}\bell_{n}^{(c,v)}(x,y,z;\lambda)=\sum_{r=0}^n\binom{n}{r}\mathcal{U}_{r}(y) \bell_{n-r}^{(c,v)}(x,z;\lambda)
		\end{equation}
		and
		\begin{equation}\label{3.4}
			{}_{\TT\!}\bell_{n}^{(s,v)}(x,y,z;\lambda)=\sum_{r=0}^n\binom{n}{r}\mathcal{U}_{r}(y) \bell_{n-r}^{(s,v)}(x,z;\lambda).
		\end{equation}
	\end{theorem}
	
	\begin{proof}
		By virtue of~\eqref{1.7}, \eqref{1.10}, and~\eqref{2.1}, we have
		\begin{equation*}
			\sum_{n=0}^\infty{}_{\TT\!}\bell_{n}^{(c,v)}(x,y,z;\lambda) \frac{t^{n}}{n!}
			=\sum_{n=0}^\infty\sum_{r=0}^n\binom{n}{r}\mathcal{U}_{r}(y) \bell_{n-r}^{(c,v)}(x,z;\lambda) \frac{t^{n}}{n!}.
		\end{equation*}
		Comparing the coefficients of $\frac{t^n}{n!}$ on both sides of the above equation, we confirm the assertion in the equation~\eqref{3.3}.
		Similarly, we can verify the assertion in the equation~\eqref{3.4}. The proof of Theorem~\ref{thm3.3-3.4} is complete.
	\end{proof}
	
	\begin{theorem}
		The parametric kinds of generalized Apostol--Bernoulli polynomials ${}_{\TT\!}\bell_{n}^{(c,v)}(x,y,z;\lambda)$ satisfy the implicit summation formula
		\begin{equation}\label{3.5}
			{}_{\TT\!}\bell_{n}^{(c,v)}(x+k,y,z;\lambda)=\sum_{r=0}^n\binom{n}{r} {}_{\TT\!}\bell_{n-r}^{(c,v)}(x,y,z;\lambda) k^{r}.
		\end{equation}
	\end{theorem}
	
	\begin{proof}
		In the equation~\eqref{2.1}, replacing $x$ by $x+k$ and using the equation~\eqref{2.1} and the series expansion of $\te^{kt}$ in the resultant equation, we obtain
		\begin{equation*}
			\sum_{n=0}^\infty{}_{\TT\!}\bell_{n}^{(c,v)}(x+k,y,z;\lambda) \frac{t^{n}}{n!}
			=\sum_{n=0}^\infty \sum_{r=0}^\infty{}_{\TT\!}\bell_{n}^{(c,v)
			}(x,y,z;\lambda) k^{r}\frac{t^{n+r}}{n!r!},
		\end{equation*}
		which, upon replacing $n$ by $n-r$ and then comparing the coefficients of $\frac{t^n}{n!}$ on both sides of the resultant equation yields the assertion in~\eqref{3.5}.
	\end{proof}
	
	\begin{remark}
		Setting $k=x$ in the equation~\eqref{3.5}, we derive
		\begin{equation*}
			{}_{\TT\!}\bell_{n}^{(c,v)}(2x,y,z;\lambda)=\sum_{r=0}^n\binom{n}{r}{}_{\TT\!}\bell_{n-r}^{(c,v)}(x,y,z;\lambda) x^{r}.
		\end{equation*}
	\end{remark}
	
	\begin{remark}
		Letting $k=1$ in the equation~\eqref{3.5}, we acquire
		\begin{equation*}
			{}_{\TT\!}\bell_{n}^{(c,v)}(x+1,y,z;\lambda)=\sum_{r=0}^n\binom{n}{r}{}_{\TT\!}\bell_{n-r}^{(c,v)}(x,y,z;\lambda).
		\end{equation*}
	\end{remark}
	
	\begin{theorem}
		The parametric kinds of generalized Apostol--Bernoulli polynomials ${}_{\TT\!}\bell_{n}^{(s,v)}( x,y,z;\lambda)$ satisfy
		\begin{equation*}
			{}_{\TT\!}\bell_{n}^{(s,v)}(x+k,y,z;\lambda)
			=\sum_{r=0}^n\binom{n}{r}\bell_{n-r}^{(s,v)}(x,z;\lambda) \TT_{r}(k,y).
		\end{equation*}
	\end{theorem}
	
	\begin{proof}
		In the equation~\eqref{2.1}, first replacing $x$ by $x+k$, then using the equations~\eqref{1.10} and~\eqref{1.22} in the resultant equation, we discover
		\begin{equation*}
			\sum_{n=0}^\infty{}_{\TT\!}\bell_{n}^{(s,v)}(x+k,y,z;\lambda) \frac{t^{n}}{n!}
			=\sum_{n=0}^\infty \sum_{r=0}^\infty\bell_{n}^{(s,v)}(x,z;\lambda)\TT_{r}(k,y) \frac{t^{n+r}}{n!r!},
		\end{equation*}
		which, upon replacing $n$ by $n-r$ and then comparing the coefficients of $\frac{t^n}{n!}$ on both sides of the resultant equation, implies the desired result.
	\end{proof}
	
	\begin{theorem}
		The parametric kinds of generalized Apostol--Bernoulli polynomials ${}_{\TT\!}\bell_{n}^{(c,v)}(x,y,z;\lambda)$ satisfy the relation
		\begin{equation}\label{3.7}
			{}_{\TT\!}\bell_{n}^{(c,v)}(\omega,y,z;\lambda)
			=\sum_{\ell=0}^n\sum_{m=0}^r\binom{n}{\ell}\binom{r}{m}(\omega -x)^{\ell+m}{}_{\TT\!}\bell_{n+r-\ell-m}^{(c,v)}(x,y,z;\lambda).
		\end{equation}
	\end{theorem}
	
	\begin{proof}
		Making use of the identity
		\begin{equation}
			\sum_{m=0}^\infty f(m) \frac{(x+y)^{m}}{m!}
			=\sum_{r=0}^\infty\sum_{s=0}^\infty f(s+r) \frac{x^{s}y^{r}}{s!r!},\label{3.8}
		\end{equation}
		in~\cite[p.~1375]{rz2}, replacing $t$ by $t+s$ in the generating function in~\eqref{2.1}, and utilizing the identity~\eqref{3.8}, we find
		\begin{equation}\label{3.9}
			\biggl(\frac{t+s}{\lambda\te^{t+s}-1}\biggr)^{v}\te^{x(t+s)} \mathcal{U}(y,t+s) \cos[z(t+s)]
			=\sum_{n=0}^\infty{}_{\TT\!}\bell_{n}^{(c,v)}(x,y,z;\lambda) \frac{(t+s)^{n}}{n!}.
		\end{equation}
		In~\eqref{3.9}, replacing $x$ by $\omega$ and then expanding the exponential function, we acquire
		\begin{equation*}
			\sum_{n=0}^\infty\sum_{r=0}^\infty{}_{\TT\!}\bell_{n+r}^{(c,v)}(\omega,y,z;\lambda) \frac{t^{n}s^{r}}{n!r!}=\sum_{k=0}^\infty \frac{(\omega-x)^{k}(t+s)^{k}}{k!}\sum_{n=0}^\infty\sum_{r=0}^\infty{}_{\TT\!}
			\bell_{n+r}^{(c,v)}(x,y,z;\lambda) \frac{t^{n}s^{r}}{n!r!}.
		\end{equation*}
		Using the identity~\eqref{3.8} on the right hand side of the above equation, replacing $n$ by $n-\ell$, $r$ by $r-m$ on the right hand side of the resultant equation give
		\begin{equation}\label{3.10}
			\sum_{n=0}^\infty\sum_{r=0}^\infty{}_{\TT\!}\bell_{n+r}^{(c,v)}(x,y,z;\lambda) \frac{t^{n}s^{r}%
			}{n!r!}=\sum_{n=0}^\infty\sum_{r=0}^\infty\sum_{\ell=0}^n\sum_{m=0}^r\frac{(\omega -x)^{\ell+m}{}_{\TT\!}\bell_{n+r-m-\ell}^{(c,v)}(x,y,z;\lambda) t^{n}s^{r}}{\ell!m!(n-\ell) !(k-m) !}.
		\end{equation}
		Comparing the coefficients of $t^n$ and $s^r$ on both sides of~\eqref{3.10} yields the assertion in~\eqref{3.7}.
	\end{proof}
	
	\begin{remark}
		Taking $z=0$ and replacing $\omega$ by $\omega+x$ in~\eqref{3.7} reduces to
		\begin{equation*}
			{}_{\TT\!}\bell_{n}^{(c,v)}(\omega+x,y,0;\lambda)
			=\sum_{\ell=0}^n\sum_{m=0}^r\binom{n}{\ell}\binom{r}{m}\omega^{\ell+m}{}_{\TT\!}\bell_{n+r-\ell-m}^{(c,v)}(x,y;\lambda).
		\end{equation*}
	\end{remark}
	
	\begin{theorem}
		Let $v, \alpha\in\mathbb{N}_{0}$. Then we have
		\begin{equation}\label{3.12}
			{}_{\TT\!}\bell_{n}^{(c,v+\alpha)}(x,y,z;\lambda)
			=\sum_{r=0}^n\binom{n}{r}\bell_{r}^{(v)}(\lambda) {}_{\TT\!}\bell_{n-r}^{(c,\alpha)}(x,y,z;\lambda)
		\end{equation}
		and
		\begin{equation}\label{3.13}
			{}_{\TT\!}\bell_{n}^{(s,v+\alpha)}(x,y,z;\lambda)
			=\sum_{r=0}^n\binom{n}{r}\bell_{r}^{(v)}(\lambda) {}_{\TT\!}\bell_{n-r}^{(s,\alpha)}(x,y,z;\lambda).
		\end{equation}
	\end{theorem}
	
	\begin{proof}
		Replacing $v$ by $v+\alpha$ in~\eqref{2.1}, we find
		\begin{align*}
			\sum_{n=0}^\infty{}_{\TT\!}\bell_{n}^{(c,v+\alpha)}(x,y,z;\lambda) \frac{t^{n}}{n!}
			&=\biggl(\frac{t}{\lambda\te^{t}-1}\biggr)^{v +\alpha}\te^{xt}\mathcal{U}(y,t) \cos(zt) \\
			&=\biggl(\frac{t}{\lambda\te^{t}-1}\biggr)^{v} \biggl(\frac{t}{\lambda\te^{t}-1}\biggr)^{\alpha}\te^{xt}\mathcal{U}(y,t) \cos(zt)
		\end{align*}
		which, upon using Apostol--Bernoulli numbers of order $v$ and~\eqref{2.1} and comparing the coefficients of $\frac{t^n}{n!}$ on both sides, implies~\eqref{3.12}.
		Similarly, we can prove the assertion in~\eqref{3.13}.
	\end{proof}
	
	\begin{theorem}
		Let $n\in\mathbb{N}_{0}$ and $\ti=\sqrt{-1}\,$. Then we have
		\begin{equation}\label{3.14}
			{}_{\TT\!}\bell_{n}^{(v)}(x+iz,y)
			={}_{\TT\!}\bell_{n}^{(c,v)}(x,y,z;\lambda)+\ti{}_{\TT\!}\bell_{n}^{(s,v)}(x,y,z;\lambda).
		\end{equation}
	\end{theorem}
	
	\begin{proof}
		Taking $z=0$ and replacing $x$ by $x+iz$ in~\eqref{2.1}, we conclude
		\begin{align*}
			\sum_{n=0}^\infty{}_{\TT\!}\bell_{n}^{(v)}(x+iz,y) \frac{t^{n}}{n!}
			&=\biggl(\frac{t}{\lambda\te^{t}-1}\biggr)^{v}\te^{(x+iz) t}\mathcal{U}(y,t)\\
			&=\biggl(\frac{t}{\lambda\te^{t}-1}\biggr)^{v}\te^{xt}\mathcal{U}(y,t) \cos(zt) +\ti\biggl(\frac{t}{\lambda\te^{t}-1}\biggr)^{v}\te^{xt}\mathcal{U}(y,t) \sin(zt)
		\end{align*}
		which, upon using~\eqref{2.2} on the right hand side and comparing the coefficients of $\frac{t^n}{n!}$ on both sides of the resultant equation, yields the assertion in~\eqref{3.14}.
	\end{proof}
	
	\begin{theorem}
		For $n\in\mathbb{N}_{0}$, we have the following summation formula
		\begin{equation}\label{3.15}
			\sum_{r=0}^n\binom{n}{r}\bell_{n-r}^{(
				v)}(x;\lambda) {}_{\TT\!}\bell_{r}^{(
				s,\beta)}(x,y,2z;\lambda)=2\underset{r=0}{\overset{n}{%
					\sum}}\binom{n}{r}\bell_{r}^{(s,v)}(
			x,z;\lambda) {}_{\TT\!}\bell_{r}^{(c,\beta
				)}(x,y,z;\lambda) .
		\end{equation}
	\end{theorem}
	\begin{proof}
		In view of the equations~\eqref{a-p} and~\eqref{2.2}, we have
		\begin{align*}
			\sum_{n=0}^\infty\bell_{n}^{(v)}(x;\lambda) \frac{t^{n}}{n!}\underset{r=0}{\overset{\infty}{\sum}}{}_{\TT\!} \bell_{r}^{(s,\beta)}(x,y,2z;\lambda) \frac{t^{r}}{r!}
			&=\biggl(\frac{t}{\lambda\te^{t}-1}\biggr)^{v}\te^{xt} \biggl(\frac{t}{\lambda\te^{t}-1}\biggr)^{\beta}\te^{xt}\mathcal{U}(y,t) \sin(2zt) \\
			&=2\biggl(\frac{t}{\lambda\te^{t}-1}\biggr)^{v}\te^{xt}\sin(zt) \biggl(\frac{t}{\lambda\te^{t}-1}\biggr)^{\beta}\te^{xt}\mathcal{U}(y,t) \cos(zt)
		\end{align*}
		which, upon using the equations~\eqref{1.7.7} and~\eqref{2.1} and comparing the coefficients of $\frac{t^n}{n!}$ on both sides of the resultant equation, yields the assertion in~\eqref{3.15}.
	\end{proof}
	
	\section{Partial derivative equations}\label{section4}
	In this section, we introduce the partial derivative equations and identities including the parametric kinds of generalized Apostol--Bernoulli polynomials ${}_{\TT\!}\bell_{n}^{(c,v)}(x,y,z;\lambda)$ and ${}_{\TT\!}\bell_{n}^{(s,v)}(x,y,z;\lambda)$. We will apply derivative operators and the generating functions in~\eqref{2.1} and~\eqref{2.2}.
	\par
	For $n\in\mathbb{N}_{0}$ and in view of the equations~\eqref{2.1} and~\eqref{2.2}, we derive the following results:
	\begin{align}
		\frac{\partial}{\partial x}\bigl[{}_{\TT\!}\bell_{n}^{(c,v)}(x,y,z;\lambda)\bigr]
		&=n\,{}_{\TT\!}\bell_{n-1}^{(c,v)}(x,y,z;\lambda), \notag \\
		\frac{\partial}{\partial x}\bigl[{}_{\TT\!}\bell_{n}^{(s,v)}(x,y,z;\lambda)\bigr]
		&=n\,{}_{\TT\!}\bell_{n-1}^{(s,v)}(x,y,z;\lambda),\notag \\
		\frac{\partial}{\partial z}\bigl[{}_{\TT\!}\bell_{n}^{(c,v)}(x,y,z;\lambda)\bigr]
		&=-n\,{}_{\TT\!}\bell_{n-1}^{(s,v)}(x,y,z;\lambda),\notag \\
		\frac{\partial}{\partial z}\bigl[{}_{\TT\!}\bell_{n}^{(s,v)}(x,y,z;\lambda)\bigr]
		&=n\,{}_{\TT\!}\bell_{n-1}^{(c,v)}(x,y,z;\lambda). \notag
	\end{align}
	The above equations show that we have
	\begin{equation*}
		\frac{\partial}{\partial x}\bigl[{}_{\TT\!}\bell_{n}^{(c,v)}(x,y,z;\lambda)\bigr]
		=\frac{\partial}{\partial z}\bigl[{}_{\TT\!}\bell_{n}^{(s,v)}(x,y,z;\lambda)\bigr]
	\end{equation*}
	and
	\begin{equation*}
		\frac{\partial}{\partial x}\bigl[{}_{\TT\!}\bell_{n}^{(s,v)}(x,y,z;\lambda)\bigr]
		=-\frac{\partial}{\partial z}\bigl[{}_{\TT\!}\bell_{n}^{(c,v)}(x,y,z;\lambda)\bigr].
	\end{equation*}
	
	\begin{theorem}
		Let $m,n\in\mathbb{N}$ and $n\ge m $. Then we have
		\begin{multline}\label{4.5}
			\frac{\partial^m}{\partial x^m}\bigl[{}_{\TT\!}\bell_{n}^{(s,v)}(x+\alpha,y,z+\beta ;\lambda)\bigr]\\
			=\sum_{r=0}^nm!\binom{n}{r}\binom{r}{m}\bigl[{}_{\TT\!}\bell_{r-m}^{(s,v)}( x,y,z;\lambda) C_{n-r}(\alpha,\beta) +{}_{\TT\!}\bell_{r-m}^{(c,v)}(x,y,z;\lambda) S_{n-r}(\alpha,\beta)\bigr].
		\end{multline}
	\end{theorem}
	
	\begin{proof}
		Replacing $x$ by $x+\alpha $ and $z$ by $z+\beta $ in~\eqref{2.2} and then applying the derivative operator $\frac{ \partial^{m}}{\partial x^{m}}$ to the resultant equation, we find
		\begin{align*}
			&\quad\sum_{n=0}^\infty\frac{\partial^{m}}{\partial x^{m}}\bigl[{}_{\TT\!}\bell_{n}^{(s,v)}(x+\alpha,y,z+\beta;\lambda)\bigr]\frac{t^{n}}{n!}\\
			&=t^{m}\biggl(\frac{t}{\lambda\te^{t}-1}\biggr)^{v}\te^{(x+\alpha) t}\mathcal{U}(y,t) \sin[(z+\beta)t]\\
			&=t^{m}\biggl[\biggl(\frac{t}{\lambda\te^{t}-1}\biggr)^{v}\te^{xt}\mathcal{U}(y,t) \sin(zt) \te^{\alpha t}\cos(\beta t) +\biggl(\frac{t}{\lambda\te^{t}-1}\biggr)^{v}\te^{xt}\mathcal{U}(y,t) \cos(zt) \te^{\alpha t}\sin(\beta t)\biggr].
		\end{align*}
		Next, in view of the equations~\eqref{1.5}, \eqref{1.6}, \eqref{2.1}, and~\eqref{2.2}, the above equation becomes
		\begin{align*}
			&\quad\sum_{n=0}^\infty\frac{\partial^{m}}{\partial x^{m}} \bigl[{}_{\TT\!}\bell_{n}^{(s,v)}(x+\alpha,y,z+\beta;\lambda)\bigr]\frac{t^{n}}{n!}\\
			&=\sum_{n=0}^\infty\sum_{r=0}^nm!\binom{n}{r}\binom{r}{m} {}{}_{\TT\!}\bell_{r-m}^{(s,v)}(x,y,z;\lambda) C_{n-r}(\alpha,\beta)
			\frac{t^{n}}{n!} \\
			&\quad+\sum_{n=0}^\infty\sum_{r=0}^nm!\binom{n}{r}\binom{r}{m}{}_{\TT\!}\bell_{r-m}^{(c,v)}(x,y,z;\lambda) S_{n-r}(\alpha,\beta) \frac{t^{n}}{n!}.
		\end{align*}
		Comparing the coefficients of $\frac{t^n}{n!}$ on both sides of the
		resultant equation yields the assertion in~\eqref{4.5}.
	\end{proof}
	
	\begin{theorem}
		Let $m,n\in\mathbb{N}$ and $ n\ge m$. Then parametric kinds of generalized Apostol--Bernoulli polynomials satisfy
		\begin{equation}\label{4.6}
			\frac{\partial^{m}}{\partial x^{m}}\bigl[{}_{\TT\!}\bell_{n}^{(c,v)}(x,y,z;\lambda)\bigr]
			=\sum_{r=0}^{n-m}m!\binom{n}{m}\binom{n-m}{r}\bell_{r}^{(\delta)}\left(\lambda\right) {}_{\TT}\!\bell_{n-r-m}^{(c,v-\delta)}(x,y,z;\lambda)
		\end{equation}
		and
		\begin{equation}\label{4.7}
			\frac{\partial^{m}}{\partial x^{m}}\bigl[{}_{\TT\!}\bell_{n}^{(s,v)}(x,y,z;\lambda)\bigr] =\sum_{r=0}^{n-m}m!\binom{n}{m}\binom{n-m}{r}\bell_{r}^{(\delta)} \left(\lambda\right){}_{\TT}\!\bell_{n-r-m}^{(s,v-\delta)}(x,y,z;\lambda).
		\end{equation}
	\end{theorem}
	
	\begin{proof}
		Using the derivative operator $\frac{\partial^{m}}{\partial x^{m}}$ in~\eqref{2.1}, we arrive at
		\begin{align*}
			\sum_{n=0}^\infty\frac{\partial^{m}}{\partial x^{m}}\biggl[{}_{\TT\!}\bell_{n}^{(c,v)}(x,y,z;\lambda)\biggr]
			&=\frac{\partial^{m}}{\partial x^{m}} \biggl[\biggl(\frac{t}{\lambda\te^{t}-1}\biggr)^{v}\te^{xt}\mathcal{U}(y,t) \cos(zt)\biggr] \\
			&=t^{m}\biggl[\biggl(\frac{t}{\lambda\te^{t}-1}\biggr)^{v}\te^{xt}\mathcal{U}(y,t) \cos(zt)\biggr] \\
			&=t^{m}\biggl[\biggl(\frac{t}{ (\lambda\te^{t}-1)}\biggr)^{\delta} \biggl(\frac{t}{\lambda\te^{t}-1}\biggr)^{v-\delta}\te^{xt}\mathcal{U}(y,t)\cos(zt)\biggr] \\
			&=\sum_{n=0}^\infty\sum_{r=0}^{n-m}m!\binom{n}{m}\binom{n-m}{r} \bell_{r}^{(\delta)}\left(\lambda\right) {}_{\TT\!}\bell_{n-r-m}^{(c,v-\delta)}(x,y,z;\lambda) \frac{t^{n}}{n!}.
		\end{align*}
		Comparing the coefficients of $\frac{t^{n}}{n!}$ on both sides, we derive the equation~\eqref{4.6}.
		Similarly, we can conclude the equation~\eqref{4.7}.
	\end{proof}

	\section{Examples and special cases}\label{section5}
	
	In this section, we demonstrate several families of polynomials, including the Apostol--Bernoulli type polynomial families, by substituting different functions for $\mathcal{U}(y,t)$, and investigate their properties.
	\par
	To give examples of the parametric kinds of generalized Apostol--Bernoulli polynomials ${}_{\TT\!}\bell_{n}^{(c,v)}(x,y,z;\lambda)$ and ${}_{\TT\!}\bell_{n}^{(s,v)}(x,y,z;\lambda)$, we recall the generating functions of some special polynomials.
	\begin{enumerate}
		\item
		The Gould--Hopper polynomials $\mathcal{H}_{n}^{(m)}(x,y)$ are defined~\cite[pp.~51--63]{t2} by
		$$
		\te^{xt+yt^m}=\sum_{n=0}^\infty\mathcal{H}_{n}^{(m)}(x,y) \frac{t^{n}}{n!},
		$$
		where $m\in\mathbb{N}$. If taking $m=2$ in the above generating function, we derive the generating function of the Hermite--Kamp\'e de Feri\'et polynomials $\mathcal{H}_{n}(x,y)$, see~\cite{t1}.
		
		\item
		The Hermite--Appell polynomials ${}_\mathcal{H}\mathcal{A}_{n}^{(m)}(x,y)$ are generated~\cite{t3} by
		$$
		\mathcal{A}(t)\te^{xt+yt^2}=\sum_{n=0}^\infty{}_\mathcal{H}\mathcal{A}_{n}^{(m)}(x,y) \frac{t^{n}}{n!}.
		$$
		What is the symbol $\mathcal{A}(t)$?
		
		\item
		The truncated exponential polynomials $e_{n}^{(m)}(x,y)$ of order $m$ are generated in~\cite{t5} by
		$$
		\te^{xt}\frac{1}{1-yt^m}=\sum_{n=0}^\infty e_{n}^{(m)}(x,y) \frac{t^{n}}{n!}.
		$$
	\end{enumerate}
	
	\subsection{Gould--Hopper--Apostol--Bernoulli type polynomials}
	When taking $\mathcal{U}(y,t)=\te^{yt^m}$ in~\eqref{2.1} and~\eqref{2.2}, we deduce the Gould--Hopper--Apostol--Bernoulli type polynomials ${}_{\mathcal{H}_{n}^{(m)}\!\!}\bell_{n}^{(c,v)}(x,y,z;\lambda)$ and ${}_{\mathcal{H}_{n}^{(m)}\!\!}\bell_{n}^{(s,v)}(x,y,z;\lambda)$ which are generated by
	\begin{equation}\label{5.1}
		\biggl(\frac{t}{\lambda\te^{t}-1}\biggr)^{v}\te^{xt+yt^m} \cos(zt)
		=\sum_{n=0}^\infty{}_{\mathcal{H}_{n}^{(m)}\!\!}\bell_{n}^{(c,v)}(x,y,z;\lambda) \frac{t^{n}}{n!}
	\end{equation}
	and
	\begin{equation}\label{5.2}
		\biggl(\frac{t}{\lambda\te^{t}-1}\biggr)^{v}\te^{xt+yt^m}\sin(zt)
		=\sum_{n=0}^\infty{}_{\mathcal{H}_{n}^{(m)}\!\!} \bell_{n}^{(s,v)}(x,y,z;\lambda) \frac{t^{n}}{n!}.
	\end{equation}
	\par
	According to~\eqref{2.3}, \eqref{2.4}, \eqref{2.5}, and~\eqref{2.6}, we see that the polynomials ${}_{\mathcal{H}_{n}^{(m)}\!\!}\bell_{n}^{(c,v)}(x,y,z;\lambda)$ and ${}_{\mathcal{H}_{n}^{(m)}\!\!}\bell_{n}^{(s,v)}(x,y,z;\lambda)$ are quasi-monomial with respect to the following multiplication and derivative operators:
	\begin{align}
		\hat{M}_{_\mathcal{H}\bell_c} &=x+ym{D_x}^{m-1}+\frac{v}{D_{x}}\frac{\lambda\te^{D_{x}}(1-D_{x})-1} {\lambda\te^{D_{x}}-1}-z\tan(zD_{x}),\label{5.3} \\
		\hat{P}_{_\mathcal{H}\bell_c} &=D_{x},\notag\\
		\hat{M}_{_\mathcal{H}\bell_s} &=x+ym{D_x}^{m-1}+\frac{v}{D_{x}}\frac{\lambda\te^{D_{x}}(1-D_{x}) -1}{\lambda\te^{D_{x}}-1}+z\cot (zD_{x}),\notag\\
		\hat{P}_{_\mathcal{H}\bell_s} &=D_{x},\notag
	\end{align}
	respectively. From~\eqref{2.8} and~\eqref{2.9}, using the derivative and multiplicative operators given above, we conclude
	\begin{equation*}
		\biggl[xD_x+ym{D_x}^{m}+v\frac{\lambda\te^{D_{x}}(1-D_{x})-1}{\lambda\te^{D_{x}}-1}-z\tan(zD_{x})D_{x} -n\biggr] {}_{\mathcal{H}_{n}^{(m)}\!\!}\bell_{n}^{(c,v)}(x,y,z;\lambda)
		=0
	\end{equation*}
	and
	\begin{equation} \label{5.8}
		\biggl[xD_x+ym{D_x}^{m}+v\frac{\lambda\te^{D_{x}}(1-D_{x}) -1}{\lambda\te^{D_{x}}-1}+z\cot(zD_{x}) D_{x}-n\biggr] {}_{\mathcal{H}_{n}^{(m)}\!\!}\bell_{n}^{(s,v)}(x,y,z;\lambda)
		=0.
	\end{equation}
	\par
	If letting $m=2$ in~\eqref{5.1} and~\eqref{5.2}, we can derive the Hermite--Kamp\'e de F\'eriet--Apostol--Bernoulli type polynomials ${}_{\mathcal{H}_{n}\hskip-3pt}\bell_{n}^{(c,v)}(x,y,z;\lambda)$ and ${}_{\mathcal{H}_{n}\hskip-3pt}\bell_{n}^{(s,v)}(x,y,z;\lambda)$. These polynomials satisfy those properties stated between~\eqref{5.3} and~\eqref{5.8}.
	
	\subsection{Hermite--Appell--Apostol--Bernoulli type polynomials}
	
	If taking $\mathcal{U}(y,t)=\mathcal{A}(t)\te^{yt^2}$ in~\eqref{2.1} and~\eqref{2.2}, then we derive the Hermite--Appell--Apostol--Bernoulli type polynomials ${}_{{}_{\mathcal{H}}{\mathcal{A}}\!}\bell_{n}^{(c,v)}(x,y,z;\lambda)$ and ${}_{{}_{\mathcal{H}}{\mathcal{A}}\!}\bell_{n}^{(s,v)}(x,y,z;\lambda)$ which are generated by
	\begin{equation*}
		\biggl(\frac{t}{\lambda\te^{t}-1}\biggr)^{v}\te^{xt+yt^2}\mathcal{A}(t) \cos(zt) =\sum_{n=0}^\infty{}_{{}_{\mathcal{H}}{\mathcal{A}}\!}\bell_{n}^{(c,v)}(x,y,z;\lambda) \frac{t^{n}}{n!}
	\end{equation*}
	and
	\begin{equation*}
		\biggl(\frac{t}{\lambda\te^{t}-1}\biggr)^{v}\te^{xt+yt^2}\mathcal{A}(t) \sin(zt)
		=\sum_{n=0}^\infty{}_{{}_{\mathcal{H}}{\mathcal{A}}\!} \bell_{n}^{(s,v)}(x,y,z;\lambda) \frac{t^{n}}{n!}.
	\end{equation*}
	\par
	From~\eqref{2.3}, \eqref{2.4}, \eqref{2.5}, and~\eqref{2.6}, we see that the polynomials ${}_{{}_{\mathcal{H}}{\mathcal{A}}\!}\bell_{n}^{(c,v)}(x,y,z;\lambda)$ and ${}_{{}_{\mathcal{H}}{\mathcal{A}}\!}\bell_{n}^{(s,v)}(x,y,z;\lambda)$ are quasi-monomial with respect to the following multiplication and derivative operators:
	\begin{align*}
		\hat{M}_{{}_\mathcal{A}\bell_c} &=x+2y{D_x}+\frac{\mathcal{A}'(D_x)}{\mathcal{A}(D_x)}+\frac{v}{D_{x}} \frac{\lambda\te^{D_{x}}(1-D_{x})-1}{\lambda\te^{D_{x}}-1}-z\tan(zD_{x}),\\
		\hat{P}_{{}_\mathcal{A}\bell_c} &=D_{x},\\
		\hat{M}_{{}_\mathcal{A}\bell_s} &=x+2y{D_x}+\frac{\mathcal{A}'(D_x)}{\mathcal{A}(D_x)}+\frac{v}{D_{x}} \frac{\lambda\te^{D_{x}}(1-D_{x})-1}{\lambda\te^{D_{x}}-1}+z\cot(zD_{x}),\\
		\hat{P}_{{}_\mathcal{A}\bell_s} &=D_{x},
	\end{align*}
	respectively.
	\par
	According to~\eqref{2.8} and~\eqref{2.9}, using the derivative and multiplicative operators given above, we obtain the differential equations
	\begin{equation*}
		\biggl[xD_x+y{D_x}^{2}+\frac{\mathcal{A}'(D_x)}{\mathcal{A}(D_x)}+v\frac{\lambda\te^{D_{x}}(1-D_{x})-1}
		{\lambda\te^{D_{x}}-1}-z\tan(zD_{x}) D_{x}-n\biggr] {}_{{}_{\mathcal{H}}{\mathcal{A}}\!}\bell_{n}^{(c,v)}(x,y,z;\lambda)
		=0
	\end{equation*}
	and
	\begin{equation*}
		\biggl[xD_x+y{D_x}^{2}+\frac{\mathcal{A}'(D_x)}{\mathcal{A}(D_x)}+v\frac{\lambda\te^{D_{x}}(1-D_{x})-1}
		{\lambda\te^{D_{x}}-1}+z\cot(zD_{x}) D_{x}-n\biggr] {}_{{}_{\mathcal{H}}{\mathcal{A}}\!}\bell_{n}^{(c,v)}(x,y,z;\lambda)
		=0.
	\end{equation*}
	\subsection{Truncated exponential Apostol--Bernoulli type polynomials}
	If $\mathcal{U}(y,t)=\frac{1}{1-yt^m}$ in~\eqref{2.1} and~\eqref{2.2}, we derive the truncated exponential Apostol--Bernoulli type polynomials ${}_{\te^{(m)}}\!\bell_{n}^{(c,v)}(x,y,z;\lambda)$ and ${}_{\te^{(m)}}\!\bell_{n}^{( s,v)}(x,y,z;\lambda)$ which are generated by
	\begin{equation*}
		\biggl(\frac{t}{\lambda\te^{t}-1}\biggr)^{v}\te^{xt}\frac{1}{1-yt^m} \cos(zt)
		=\sum_{n=0}^\infty{}_{\te^{(m)}}\!\bell_{n}^{(c,v)}(x,y,z;\lambda)\frac{t^{n}}{n!}
	\end{equation*}
	and
	\begin{equation*}
		\biggl(\frac{t}{\lambda\te^{t}-1}\biggr)^{v}\te^{xt}\frac{1}{1-yt^m} \sin(zt)
		=\sum_{n=0}^\infty{}_{\te^{(m)}}\!\bell_{n}^{( s,v)}(x,y,z;\lambda)\frac{t^{n}}{n!}.
	\end{equation*}
	\par
	Considering~\eqref{2.3}, \eqref{2.4}, \eqref{2.5}, and~\eqref{2.6}, we see easily that the polynomials ${}_{\te^{(m)}}\!\bell_{n}^{(c,v)}(x,y,z;\lambda)$ and ${}_{\te^{(m)}}\!\bell_{n}^{(s,v)}(x,y,z;\lambda)$ are quasi-monomial with respect to the following multiplication and derivative operators:
	\begin{align*}
		\hat{M}_{{}_{\te^{(m)}}\!\bell_c} &=x+\frac{v}{D_{x}}\frac{\lambda
			\te^{D_{x}}(1-D_{x})-1}{\lambda\te^{D_{x}}-1}+\frac{my{D_{x}^{m-1}}}{1-yD_{x}^m}-z\tan (zD_{x}),\\
		\hat{P}_{{}_{\te^{(m)}}\!\bell_c} &=D_{x},\\
		\hat{M}_{{}_{\te^{(m)}}\!\bell_s} &=x+\frac{v}{D_{x}}\frac{\lambda
			\te^{D_{x}}(1-D_{x}) -1}{\lambda\te^{D_{x}}-1}+\frac{my{D_{x}^{m-1}}}{1-yD_{x}^m}+z\cot (zD_{x}),\\
		\hat{P}_{{}_{\te^{(m)}}\!\bell_s} &=D_{x},
	\end{align*}
	respectively.
	\par
	From~\eqref{2.8} and~\eqref{2.9}, considering the derivative and multiplicative operators given above, we acquire
	\begin{equation*}
		\biggl[xD_x+v\frac{\lambda\te^{D_{x}}(1-D_{x}) -1}{\lambda\te^{D_{x}}-1}+\frac{my{D_{x}^{m}}}{1-yD_{x}^m}-z\tan(zD_{x}) D_{x}-n\biggr] {}_{\te^{(m)}}\!\bell_{n}^{(c,v)}(x,y,z;\lambda)
		=0
	\end{equation*}
	and
	\begin{equation*}
		\biggl[xD_x+v\frac{\lambda\te^{D_{x}}(1-D_{x}) -1}{\lambda\te^{D_{x}}-1}+\frac{my{D_{x}^{m}}}{1-yD_{x}^m}+z\cot(zD_{x}) D_{x}-n\biggr] {}_{\te^{(m)}}\!\bell_{n}^{(s,v)}(x,y,z;\lambda)
		=0.
	\end{equation*}
	
	\section{Conclusion}
	In this paper, we have introduced the new generating functions using Euler's formula  for parametric kinds of generalized Apostol--Bernoulli polynomials which are of quasi-monomial properties. Then we have investigated quasi-monomial properties and identities, some addition and summation formulas and partial derivative formulas. Then, in the special cases of $\mathcal{U}\left(y,t\right)$, we construct new subpolynomial families such as the Gould--Hopper--Apostol--Bernoulli type polynomials, the Hermite--Appell--Apostol--Bernoulli type polynomials and truncated exponential Apostol--Bernoulli type polynomials. Lastly, we investigated the multiplicative and derivative operators and differential equations of these new families.
	\par
	This paper is a slightly revised version of the preprint~\cite{apostol-bernoulli-polynomials.tex}.

	\end{document}